\documentclass[10pt]{article}

\usepackage{amsmath,amssymb,amscd}
\usepackage{cite}
\usepackage{graphicx,epstopdf}
\usepackage{float}
\usepackage{authblk} 
\usepackage{epigraph}
\numberwithin{equation}{section}
\usepackage[usenames]{color}
\usepackage{color}
\usepackage{colortbl}
\usepackage[linktocpage=true,plainpages=false,pdfpagelabels=false]{hyperref}
\hypersetup{colorlinks=true, urlcolor=blue, citecolor=blue, linkcolor=blue}

\usepackage{soul}
\usepackage{tikz-cd}

\def\e{\epsilon}

\def\O{\Omega}

\def\bp{\begin{proposition}}
\def\ep{\end{proposition}}
\def\bt{\begin{theo}}
\def\et{\end{theo}}
\def\be{\begin{equation}}
\def\ee{\end{equation}}
\def\bl{\begin{lemma}}
\def\el{\end{lemma}}
\def\bc{\begin{corollary}}
\def\ec{\end{corollary}}

\def\bd{\begin{definition}}
\def\ed{\end{definition}}

\def\max{{\rm max\,}}

\def\max{{\rm max\,}}

\newcommand{\tab}{\hspace*{2em}}
\newtheorem{theo}{Theorem}[section]
\newtheorem{lemma}{Lemma}[section]
\newtheorem{definition}{Definition}[section]
\newtheorem{corollary}{Corollary}[section]
\newtheorem{proposition}{Proposition}[section]

\begin{document}

\title{Accuracy of noisy Spike-Train Reconstruction: a Singularity Theory point of view}

\author[1]{Gil Goldman}
\author[2]{Yehonatan Salman}
\author[3]{Yosef Yomdin}
\date{} 

\affil[1,2,3]{Department of Mathematics, The Weizmann Institute of Science, Israel}

\medskip

\affil[1]{Email: gilgoldm@gmail.com}
\affil[2]{Email: salman.yehonatan@gmail.com}
\affil[3]{Email: yosef.yomdin@weizmann.ac.il}

\maketitle

{\it To the Memory of Egbert Brieskorn.
Among the most important events which inspired scientific interests of the third author was (in late 1960s) Brieskorn's discovery that 28 Milnor's exotic spheres can be
described by simple algebraic equations, and (in early 1980s) participation in Brieskorn Singularities seminar in Bonn.}

\begin{abstract}

This is a survey paper discussing one specific (and classical) system of algebraic equations - the so called ``Prony system''. We provide a short overview of its unusually wide connections with many different fields of Mathematics, stressing the role of Singularity Theory. We reformulate Prony System as the problem of reconstruction of ``Spike-train'' signals of the form $F(x)=\sum_{j=1}^d a_j\delta(x-x_j)$ from the noisy moment measurements. We provide an overview of some recent results of \cite{akinshin2015accuracy,akinshin2017accuracy,akinshin2017error,goldman2018prony,batenkov2017accurate,batenkov2016stability,batenkov2013accuracy,batenkov2013geometry,yomdin2010singularities} on the ``geometry of the error amplification'' in the reconstruction process, in situations where the nodes $x_j$ near-collide. Some algebraic-geometric structures, underlying the error amplification, are described (Prony, Vieta, and Hankel mappings, Prony varieties), as well as their connection with Vandermonde mappings and varieties. Our main goal is to present some promising fields of possible applications of Singulary Theory.

\end{abstract}



\section{Introduction}\label{Sec:Intro}
\setcounter{equation}{0}

In this paper we consider the
classical {\it Prony system} of algebraic equations, with the real unknowns $a_j,x_j, \ j=1,\ldots,d,$ and with
the right hand side formed by the known real ``measurements'' $m_0,\ldots,m_{2d-1}$. This system has a form

\be\label{eq:Prony.system1}
\sum_{j=1}^d a_j x_j^k = m_k, \ k= 0,1,\ldots,2d-1.
\ee
We denote by $A=(a_1,\ldots,a_d) \in {\mathbb R}^d$ and $X=(x_1,\ldots,x_d) \in {\mathbb R}^d, \ x_1\le x_2\le \ldots \le x_d,$
the unknowns in system (\ref{eq:Prony.system1}), and denote by ${\cal P}^A_d$ (resp. ${\cal P}^X_d$) the ``parameter spaces'' of the unknowns $A$ and $X,$ respectively. ${\cal P}_d={\cal P}^A_d\times {\cal P}^X_d$ denotes the total parameter space of $(A,X)$.
The space (isomorphic to ${\mathbb R}^{2d}$) of the right-hand sides $\mu=(m_0,m_1,\ldots,m_{2d-1})$ of (\ref{eq:Prony.system1}) is denoted by ${\cal M}_d$.

\smallskip

In what follows we will usually identify $(A,X)$ with a ``spike-train signal''

\be\label{eq:s.p.signal}
	F(x)=\sum_{j=1}^{d}a_{j}\delta(x-x_{j}).
\ee
Clearly, the moments $m_k(F)=\int x^k F(x)dx, \ k=0,1,\ldots,$ are given by $m_k(F)=\sum_{j=1}^d a_j x_j^k$, so reconstructing $F$ from its $2d-1$ initial moments is equivalent to solving (\ref{eq:Prony.system1}), with $m_k=m_k(F)$.



\medskip


Prony system appears in many classical theoretical and applied mathematical problems. In Section \ref{Sec:Appear.Prony} we discuss some of these appearances. Explicit solution of
(\ref{eq:Prony.system1}) was given already in \cite{prony1795essai} (see Section \ref{Sec:Exp.Sol} below).

There exists a vast literature on Prony and similar systems. In particular, the bibliography in \cite{auton1981investigation} (1981) contains more than 50 pages. Most of recent applications are in Signal Processing. As a very partial sample we mention that in \cite{blu2008sparse} and in many other publications a method, essentially equivalent to solving Prony system, was used in reconstructing signals with a ``finite rate of innovation''. In \cite{peter2013generalized,plonka2014prony} the applicability of Prony-type systems was extended to some new wide and important classes of signals.
In \cite{bernardi2014comparison,comon2008symmetric} multidimensional Prony systems were investigated via symmetric tensors, in particular, connecting them to the polynomial Waring problem. In \cite{eldar2015sampling} Prony system appears in a general context of Compressed Sensing. In \cite{batenkov2012algebraic,batenkov2015complete} Prony-like systems were used in reconstructing piecewise-smooth functions from their Fourier data. Finally, in \cite{batenkov2015complete} the same reconstruction accuracy as for smooth functions was demonstrated (thus confirming the Eckhoff conjecture).

\smallskip

Some applications of Prony system are of major practical importance, and various algorithms and numerical methods have been developed for its solution (see \cite{potts2015Fast} and references therein). However, in a (very important) case when some of the nodes $x_j$ nearly collide, while the measurements are noisy, these collision singularities lead to major mathematical and numerical difficulties. In particular, this happens in the context of the ``super-resolution problem'', which was investigated in many recent publications. See \cite{akinshin2015accuracy,akinshin2017accuracy,akinshin2017error,batenkov2017accurate,batenkov2016stability,batenkov2013accuracy,batenkov2013geometry,candes2013super,candes2014towards,demanet2015recoverability,donoho1992superresolution,fernandez2016super,morgenshtern2016super} as a small sample.

\smallskip

Notice that the Prony system \eqref{eq:Prony.system1} is linear (with the Vandermonde matrix on the ``nodes'' $X$) with respect to the ``amplitudes'' $A$, while it is highly
nonlinear with respect to $X$.
As the nodes collide (or near-collide), the Vandermonde determinant vanishes. Even knowing the position of the nodes, the reconstruction of the amplitudes is still ill-posed.

Thus singularities enter the solution process of the Prony system because of its geometric nature, no matter what solution method do we use. We believe that using the tools of Singularity Theory in this problem is well justified. In \cite{yomdin2010singularities,batenkov2013geometry} we study the algebraic nature of nodes collision. In particular, we include into consideration the ``confluent Prony systems'', corresponding to signals with multiple nodes, and with the derivatives of the $\delta$-function. We also introduce and study in \cite{batenkov2013geometry} the ``bases of finite differences'' in the signal space ${\cal P}_d$, which behave coherently as the nodes collide.

\smallskip

In the present paper we give, following \cite{akinshin2015accuracy,akinshin2017accuracy,akinshin2017error,goldman2018prony,batenkov2017accurate,batenkov2016stability,batenkov2013accuracy,batenkov2013geometry,yomdin2010singularities}, a somewhat different point of view on the problem, stressing the role of Singularity Theory in understanding of Prony systems {\it with noisy right-hand side}. Below we discuss the following main topics:

\smallskip

\noindent 1. In case of near-colliding nodes the initial measurements errors may be strongly amplified in the solution, making it unfeasible. However, the possible error-affected solutions are not distributed uniformly, but rather tightly concentrated along certain algebraic sets, known a priori (``Prony varieties'' - see Sections \ref{Sec:Error.Amplif.Prony} and \ref{Sec:Prony.Var.in.Three.Spaces} below).

{\it Prony varieties are generalizations (via ``making free'' the amplitudes $A$) of the Vandermonde mappings and varieties, introduced and studied in Singularity Theory in \cite{arnol1986hyperbolic,kostov1989geometric} and other publications (see Section \ref{Sec:Pro.Vie.Han.mappings}}).

\smallskip

\noindent 2. A related notion is that of ``Prony scenarios'' (Section \ref{Sec:Prony.Scenarios}), which predict the error behavior along the Prony curves. {\it An important part here is the description of the combinatorics of real zeroes in polynomial pencils, actively studied in Singularity Theory - see  \cite{borcea2004classifying,kostov2002root,kostov2011topics,kurdyka201nuij}}.

\smallskip

\noindent 3. In the presence of the measurement noise, statistical estimations of feasible solutions can be used. These methods are considered in the literature as superior in accuracy, but their practical implementation is difficult, because of complicated nonlinear minimization problems involved. We expect that the tools developed in Singularity Theory for the study of ``maxima of smooth functions'', ``cut-loci'', and similar objects, can be useful here (Section \ref{Prob.setting}).

\smallskip

\noindent 4. In the case of the real nodes $X$ (mainly presented in this paper) hyperbolic polynomials become a central topic in all the problems above. {\it Hyperbolic polynomials and related objects are actively studied in Singularity Theory (see \cite{arnol1986hyperbolic,borcea2004classifying,dimitrov2010distances,kostov1989geometric,kostov2002root,kostov2011topics,kurdyka201nuij} as a partial sample), and we expect some of the available results to be directly applicable to Prony systems}.

\medskip

Among other common topics with Singularity Theory we shortly discuss below rank stratification of the space of Hankel-type matrices, solving parametric linear systems, polynomial Waring problem, and finite differences. We hope that the connections presented will proof useful in both domains.

\section{Some appearances of the Prony system}\label{Sec:Appear.Prony}
\setcounter{equation}{0}

We outline here some prominent classical appearances of the Prony system.

\subsection {Exponential Interpolation}\label{Sec:Exp.Interpol}

This was the problem studied by Prony himself in \cite{prony1795essai}. We consider an interpolation problem for a given function $f(x)$ at the $2d$ consequent integer points $0,1,\ldots,2d-1$, with the interpolant being the sum of the exponents

$$
\sum_{j=1}^d a_je^{\zeta_j x}.
$$

We can choose freely $2d$ parameters $a_j,\zeta_j$, in order to fit the values $y_k=f(k), \ k=0,\ldots,2d-1.$ Substituting $x=k$, and denoting $e^{\zeta_j}$ by $x_j$ we get the Prony system of equations

$$
\sum_{j=1}^d a_je^{k\zeta_j}=\sum_{j=1}^d a_jx_j^k=y_k, \ k=0,1,\ldots,2d-1.
$$

\subsection{Gauss quadratures}\label{Sec:Gauss quadratures}

Let $\lambda$ be a measure on the real line $\mathbb R$. For a given $d$ we want to find $d$ points $x_1,\ldots,x_d \in {\mathbb R}$, and $d$ real coefficients $a_1,\ldots,a_d$ such that the quadrature formula

\be\label{eq:Gauss}
\int g(x) d\lambda \approx \sum_{j=1}^d a_jg(x_j)
\ee
be accurate for $g$ being any polynomial of degree at most $2d-1$. By linearity, it is sufficient to get an equality in (\ref{eq:Gauss}) only for $g$ being the monomials $x^k, \ k=0,1,\ldots,2d-1$, and this leads immediately to the Prony system

\be\label{eq:Gauss.Prony}
\sum_{j=1}^d a_jx_j^k=m_k(\lambda): =\int x^k d\lambda, \ k=0,1,\ldots,2d-1,
\ee
with the right-hand side given by the consecutive moments $m_k(\lambda)$ of the measure $\lambda$.

\smallskip

Another interpretation is that we are looking for an atomic measure (a spike-train signal) $\tilde \lambda = \sum_{j=1}^{d}a_{j}\delta(x-x_{j})$ satisfying $m_k(\tilde \lambda)=m_k(\lambda), \ k=0,1,\ldots,2d-1.$


\subsection{Moment Theory and Pad\'e approximations}\label{Sec:Moments.Pade}

The classical Hamburger Moment problem consists in providing necessary and sufficient conditions for a sequence $m=\{m_0,m_1,\ldots,m_k,\ldots\}$ to be the sequence of the consecutive moments $m_k=m_k(\lambda)=\int x^k d\lambda, \ k=0,1,\ldots,$ of a non-atomic positive measure $\lambda$ on the real line $\mathbb R$, and in reconstructing $\lambda$ from $m$. The condition is that all the Hankel-type matrices

\begin{equation}\label{eq:Hankel.Matrix.1}
M_d(m)=				
\begin{bmatrix}
m_{0}  & m_1 &... & m_{d-1}     \\
m_{1}  & m_2 &... & m_{d}       \\
\reflectbox{$\ddots$}  & \reflectbox{$\ddots$} &     &        \\					
m_{d-1}& m_d   & ...    & m_{2d-2}
\end{bmatrix}
\tab d=0,1,...,	
\end{equation}	
are positive definite. The proof essentially consists in Gaussian quadrature approximation of the measure $\lambda$ by positive atomic measures $\lambda_d = \sum_{j=1}^{d}a_{d,j}\delta(x-x_{d,j}), \ d=0,1,\ldots$, satisfying the condition $m_k(\lambda_d)=m_k, \ k=0,1,\ldots,2d-1,$ i.e. solving the Prony systems

\be\label{eq:Moments.Prony}
\sum_{j=1}^d a_{d,j}x_{d,j}^k=m_k, \ k=0,1,\ldots,2d-1, \ d=0,1,\ldots,
\ee
with the right-hand side given by the input sequence $m=\{m_0,\allowbreak m_1,\allowbreak \ldots,m_k,\ldots\}$.

\smallskip

Another point of view is provided by the Pad\'e approximation approach. For a sequence $m$ as above consider a formal power series at infinity

\be\label{eq:formal.ser}
f(z)=\sum_{k=0}^\infty m_kz^{-k-1}.
\ee
The $d$-th (diagonal) Pad\'e approximant of $f(z)$ is a rational function $R_d(z)=\frac{P_d(z)}{Q_d(z)}$ with $P_d,Q_d$ polynomials in $z$ of the degrees $d-1$ and $d$, respectively, such that the Taylor development of $R_d(z)$ at infinity has the form

\be\label{eq:rat.ser}
R_d(z)=\sum_{k=0}^{2d-1} m_kz^{-k-1} + O(z^{-2d-1}).
\ee
In other words, the first $2d$ Taylor coefficients of $R_d(z)$ are $m_0,\ldots,m_{2d-1}.$

\smallskip

Write $R_d(z)$ as the sum of elementary fractions, and develop at infinity:

$$
R_d(z)=\sum_{j=0}^d \frac{a_{d,j}}{z-x_{d,j}}=\sum_{j=0}^d \frac{a_{d,j}}{z(1-\frac{x_{d,j}}{z})}=\sum_{j=0}^d \frac{a_{d,j}}{z}(1+\frac{x_{d,j}}{z}+(\frac{x_{d,j}}{z})^2+\ldots)=
$$
$$
= \sum_{k=0}^\infty \tilde m_kz^{-k-1},
$$
where
$$
\tilde m_k = \sum_{j=1}^d a_{d,j}x_{d,j}^k, \ k=0,1,\ldots .
$$
Thus condition (\ref{eq:rat.ser}) becomes the Prony system \eqref{eq:Moments.Prony}.

\smallskip

We do not discuss here other remarkable connections of the Prony system, provided by the classical Moment Theory, in particular,
with continued fractions and orthogonal polynomials, see, for example, \cite{nikishin1991rational}.






\subsection{Polynomial Waring problem}\label{Sec:Pol.Waring}

We consider only the case of two variables (in more variables the calculations are, essentially, the same). Let $P(x,y)=\sum_{i=0}^m b_i x^{m-i}y^i$ be a homogeneous polynomial of degree $m$ in $(x,y)$. We look for a representation of $P$ as a sum of $m$-th powers of $d$ linear forms in $(x,y)$:

\be\label{eq:Waring.setting}
P(x,y)=\sum_{j=1}^d (\eta_jx+\zeta_jy)^m,
\ee
within an attempt to minimize $d$ in this expression. This problem is actively studied today. Many important results on generic and non-generic configurations in different degrees and dimensions are available. For details we refer the reader to \cite{bernardi2014comparison,ciliberto2001geometric,comon2008symmetric,lundqvist2017generic,miranda1999linear}, and references therein, as a very partial sample.


\smallskip

Let us put $x=1$ in (\ref{eq:Waring.setting}). We get

\be\label{eq:Waring.setting1}
P(1,y)=\sum_{i=0}^m b_i y^i=\sum_{j=1}^d (\eta_j+\zeta_jy)^m=\sum_{j=1}^d \eta_j(1+\frac{\zeta_j}{\eta_j}y)^m.
\ee
Denoting in (\ref{eq:Waring.setting1}) the fraction $\frac{\zeta_j}{\eta_j}$ by $\xi_j$ we get

$$
\sum_{i=0}^m b_i y^i=\sum_{j=1}^d \eta_j(1+\xi_j y)^m=\sum_{j=1}^d \eta_j\sum_{i=0}^m (^d_i)\xi_j^iy^i=\sum_{i=0}^m y^i \sum_{j=1}^d \eta_j (^d_i)\xi_j^i.
$$
Comparing the coefficients of $y^i$ on the two sides we obtain
$$
\sum_{j=1}^d \eta_j (^d_i)\xi_j^i=b_i, \ i=0,\ldots,m.
$$
Finally, dividing by $(^d_i)$ and denoting $b_i/(^d_i)$ by $\mu_i$, we get the Prony system $\sum_{j=1}^d \eta_j\xi_j^i=\mu_i, \ i=0,\ldots,m.$

\section{Explicit solution of the Prony system}\label{Sec:Exp.Sol}
\setcounter{equation}{0}

From now on, and till Section \ref{Sec:Solv.Real.Hyperb}, we allow complex nodes and amplitudes $(A,X)$. In Section \ref{Sec:Solv.Real.Hyperb} we return to the real case, and explain the role of hyperbolic polynomials in the solution process.

\smallskip

In order to solve explicitly Prony system

\be\label{eq:Moments.Prony11}
\sum_{j=1}^d a_{j}x_{j}^k=m_k, \ k=0,1,\ldots,2d-1,
\ee
consider the $d$-th diagonal Pad\'e approximant $R_d(z)$ of the moment generating function, defined by (\ref{eq:rat.ser}) above.




Writing $R_d(z)$ as $R_d(z)=\frac{P_d(z)}{Q_d(z)}$ with
$$
P_d(z)=b_0+b_1z+\ldots + b_{d-1}z^{d-1}, \ Q_d(z)=c_0+c_1z+\ldots + c_{d-1}z^{d-1}+z^d,
$$
substituting into (\ref{eq:rat.ser}), and comparing coefficients, we obtain the following linear system of equations for the coefficients $c=(c_0,\ldots,c_{d-1})$ of the denominator $Q$:

\begin{equation}\label{eq.solving.for.Q}				
			\begin{bmatrix}
				m_{0}  & m_1 &... & m_{d-1}     \\
				m_{1}  & m_2 &... & m_{d}       \\
				\reflectbox{$\ddots$}  & \reflectbox{$\ddots$} &     &        \\					
				m_{d-1}& m_d   & ...    & m_{2d-2}
			\end{bmatrix}	
			\begin{bmatrix}
				c_0\\
				c_1\\
				\vdots\\
				c_{d-1}
			\end{bmatrix}
			=
			-
			\begin{bmatrix}
				m_d\\
				m_{d+1}\\
				\vdots\\
				m_{2d-1}
			\end{bmatrix}	.	
		\end{equation}	
with the Hankel matrix $M_d(\mu), \ \mu=(m_0,\ldots,m_{2d-1}).$

\smallskip

Finding $c$ from (\ref{eq.solving.for.Q}), we then find the coefficients $b=(b_0,\ldots,b_{d-1})$ of the numerator $P$ as
$$
b_0=m_0c_0, \ \ \ b_1=m_0c_1+m_1c_0, \ \ldots , \ \ b_{d-1}=m_0c_{d-1}+\ldots +m_dc_0.
$$
This provides us explicitly the Pad\'e approximant $R_d(z)=\frac{P(z)}{Q(z)}.$ In order to find $a_j,x_j$ it remains to represent $R_d$ as the sum of the elementary fractions $R_d(z)=\sum_{j=0}^d \frac{a_{j}}{z-x_{j}}$. Essentially, this procedure appeared already in the Prony paper \cite{prony1795essai}, and it remains a basis for most of recent algorithms.

\subsection{Solvability conditions}\label{Sec:Solv.Cond}

Solvability conditions for (\ref{eq.solving.for.Q}) (and for the Prony system) are well known in the classical Moment Theory, in Pad\'e approximations, and in other related fields, sometimes in quite different forms. One of possible formulations, convenient for our setting, was given in \cite{batenkov2013geometry}. In order to present these conditions in a compact form, we allow complex nodes and amplitudes, as well as multiple nodes. (Including multiple nodes requires a rather accurate treatment, which we omit here. Details are given in \cite{batenkov2013geometry}).


\smallskip

From the right hand side $\mu =(m_0,\ldots,m_{2d-1}) \in {\cal M}_d$ we form the extended $d\times (d+1)$ Hankel matrix $\tilde {\cal M}_d(\mu)$.

\bt\label{th:solvsb.cond}(See \cite{batenkov2013geometry}).
Prony system (\ref{eq:Moments.Prony11}) is solvable if and only if the following condition is satisfied: let the rank of $\tilde {\cal M}_d(\mu)$ be equal to $r\leq d$. Then the left-upper $r\times r$ minor of $\tilde {\cal M}_d(\mu)$ is non-zero.
\et
Thus solvability of (\ref{eq:Moments.Prony11}) can be read out from the right-hand side $\mu$ through the ``rank stratification $\Sigma$'' of the moment space ${\cal M}_d$.

\medskip

Rank stratification for various classes of matrices is very important in Singularity Theory, and an extensive literature exists on this topic. Let us mention just \cite{mather1973solutions,fukuda2004singularities}, which may be directly related to our study of Prony system. Specifically, J. Mather's theorem in \cite{mather1973solutions} provides conditions for existence of smooth (in parameters) solutions of parametric families of linear systems (see also related results in \cite{fukuda2004singularities}). We expect that Mather's theorem can be applied to the above system (\ref{eq.solving.for.Q}), providing a very important information on the behavior of solutions of (\ref{eq.solving.for.Q}) as $\mu$ approaches the low rank strata of $\Sigma$.

Let us mention also \cite{goryunov1995semi,marar1989multiple,Ballesteros2017Sing} where finite differences, and semi-simplicial resolutions, appear in study of Image singularities. They may be related to the study of the Prony mapping, via bases of finite differences in \cite{batenkov2013geometry}.



\smallskip

\section{Prony, Vieta and Hankel mappings}\label{Sec:Pro.Vie.Han.mappings}
\setcounter{equation}{0}

In this section we suggest an algebraic-geometric picture capturing, to some extent, the mathematical structure of the solution procedure in Section \ref{Sec:Exp.Sol}. An important fact is that this picture appears as a natural extension of a construction, well known in Singularity Theory: that of Vandermonde mapping and Vandermonde varieties, developed by Arnold, Givental, Kostov and others in the 1980's (see \cite{arnol1986hyperbolic,kostov1989geometric,froberg2016vandermonde} and references therein).

Consider the following mappings:

\smallskip

\noindent 1. The Prony map:
$$
PM: {\cal P}_d\to {\cal M}_d, \ PM(F)=(m_0(F),\ldots,m_{2d-1}(F)).
$$
For each fixed amplitudes $A=(a_1,\ldots,a_d)$ the restriction of the Prony map to $A\times {\cal P}_d^X$ coincides with the corresponding Vandermonde map, as defined in \cite{arnol1986hyperbolic,kostov1989geometric,froberg2016vandermonde}.

\smallskip

We call the space of all monic polynomials of degree $d$, $Q(z)=c_0+c_1z+\ldots + c_{d-1}z^{d-1}+z^d$,  the polynomial space ${\cal V}_d$.

\smallskip

\noindent 2. The Vieta map:
$$
VM: {\cal P}_d\to {\cal V}_d, \ VM(F)= Q_F(z)=z^d+\sigma_1(F)z^{d-1}+\ldots + \sigma_1(F).
$$
Here $\sigma_i(F)=\sigma_i(x_1,\ldots,x_d)$ is the $i$-th symmetric polynomial in the nodes $X$ of $F$, and $Q(z)=Q_F(z)$ is the normalized polynomials with the roots $x_1,\ldots,x_d$. Notice that the Vieta map depends only on the nodes $X$ of $F$, but not on its amplitudes $A$.


\smallskip

\noindent 3. The Hankel map:
$$
HM: {\cal M}_d\to {\cal V}_d.
$$
This map associates to any $\mu= (m_0,\ldots,m_{2d-1})\in {\cal M}_d$ the polynomial $Q\in {\cal V}_d$ obtained through solving a linear system
(\ref{eq.solving.for.Q})

\smallskip

Notice that in the coordinates $A,X$ in the signal space ${\cal P}_d$ the mappings $PM$ and $VM$ are polynomial, while the mapping  $H$ in the coordinates $\mu=(m_0,\ldots,m_{2d-1})$ is rational, with the denominator $\Delta(\mu)=\det M_d(\mu)$, as provided by the Cramer rule.

\smallskip

We can put the mappings above into a mapping diagram $D$:

\smallskip

$$
	\begin{tikzcd}[column sep=normal]
	& {\cal V}_d  & \\
	{\cal P}_d \arrow{ur}[swap]{VM} \arrow{rr}[swap]{PM} & & {\cal M}_d \arrow{ul}{HM}
	\end{tikzcd}
$$

\smallskip

Now, a simple and basic fact, expressing the Prony solution algorithm, is the following:

\bp\label{prop:commut.diagr}
The mapping diagram $D$ is commutative, i.e.
$$
VM = HM\circ PM.
$$
\ep
The proof was, essentially given in Section \ref{Sec:Exp.Sol} above.


\smallskip

The role of each of the three spaces in the solution process is different, and some important structures may look quite differently in these spaces. Below we give some examples.




\section{Prony varieties}\label{Sec:Prony.Var.in.Three.Spaces}
\setcounter{equation}{0}

In this section we define, following \cite{akinshin2015accuracy,akinshin2017accuracy,akinshin2017error}, the ``Prony varieties'', which play an important role in the description of error amplification in solving Prony system.

\smallskip

Possessing the diagram $D$ we can choose the easiest place to define the Prony varieties, which is the moment space ${\cal M}_d$. For each $d\le q \le 2d-1$ and for a given $\mu\in {\cal M}_d$, the ``moment Prony variety''
$S^{\cal M}_q(\mu)$ is the coordinate subspace in ${\cal M}_d$, passing through the point $\mu$, where the first $q+1$ moments $m_0,\ldots,m_q$ are constant.

\smallskip

The ``signal Prony variety'' $S^{\cal P}_q(\mu)$ is the preimage under the Prony mapping $PM$ of the moment Prony variety $S^{\cal M}_q(\mu)$. Thus in ${\cal P}_d$ this variety is defined by the system of equations

\be\label{eq:Moments.Prony17}
\sum_{j=1}^d a_{j}x_{j}^k=m_k, \ k=0,1,\ldots,q,
\ee
which is formed by the first $q+1$ equations of the complete Prony system (\ref{eq:Moments.Prony11}). This was the original definition of the ``Prony leaves'' in \cite{akinshin2015accuracy} and in later publications.
For each fixed amplitudes $A=(a_1,\ldots,a_d)$ the signal Prony variety, intersected with $A\times {\cal P}_d^X$,  coincides with the corresponding Vandermonde variety, as defined in \cite{arnol1986hyperbolic,kostov1989geometric,froberg2016vandermonde}.
We believe that the results of these papers may be important in study of Prony varieties, and we give more detail in \cite{goldman2018prony}.

\smallskip

The ``polynomial Prony variety'' $S^{\cal V}_q(\mu)\subset {\cal V}_d$ is the image under the Hankel map $HM$ of the moment Prony variety $S^{\cal M}_q(\mu)$. We have the following fact:

\bp\label{prop:linearity.Prony.var}(\cite{goldman2018prony})
For $q \ge d$ the polynomial Prony varieties $S^{\cal V}_q(\mu)$ are affine subspaces in ${\cal V}_d$, defined by the linear equations

\be\label{eq:linear.eq.Pr.var}
\begin{array}{c}
\mu_{d-1}c_1+\mu_{d-2}c_2+\ldots+\mu_0c_d =-\mu_d\\
\mu_{d}c_1+\mu_{d-1}c_2+\ldots+\mu_1c_d =-\mu_{d+1}\\

..........\\
\mu_{q-1}c_1+\mu_{q-2}c_2+\ldots+\mu_{q-d}c_d =-\mu_{q} .\\
\end{array}
\ee
\ep


\smallskip

In the signal space we obtain in \cite{goldman2018prony} the following description of the (node projections) of the Prony varieties $S^{\cal P}_q, \ d\le q \le 2d-1$:

\bt\label{eq:eq.prony.x}(\cite{goldman2018prony})
The projection $S^{{\cal P},X}_q(\mu)$ of the signal Prony variety $S^{\cal P}_q(\mu)$ to the nodes space ${\cal P}^X_d$ is defined in ${\cal P}^X_d$ by the equations (\ref{eq:linear.eq.Pr.var}), with $c_j, \ j=1,\ldots,d,$ replaced by the symmetric polynomials $\sigma_j(x_1,\ldots,x_d)$.

\smallskip

In the real case, the Vieta map $VM$ provides a diffeomorphism of the interior of $S^{{\cal P},X}_q(\mu)$ to the interior of the intersection of the polynomial Prony varieties $S^{\cal V}_q(\mu)$ with the set $H_d$ of hyperbolic polynomials in ${\cal V}_d$. The inverse is given by the ``root mapping'' $RM$, which associates to a hyperbolic polynomial $Q\in H_d^\circ$ its ordered roots $x_1<\ldots < x_d$.
\et

\smallskip

For any $q$ between $d$ and $2d-2$ we can consider the parametrization of the polynomial Prony varieties $S^{\cal V}_q$ through the last ``free'' moments $m_{q+1},\ldots,m_{2d-1}$ in the right hand side of (\ref{eq.solving.for.Q}). This is the restriction of the mapping $HM$ to the the moment Prony varieties $S^{\cal M}_q(\mu)$, i.e., to the coordinate subspaces in ${\cal M}_d$, passing through the point $\mu$, where the first $q$ moments $m_0,\ldots,m_q$ are constant. We have:

\bp\label{prop:degree.of.param}(\cite{goldman2018prony})
The restriction of the mapping $HM$ to the the moment Prony varieties $S^{\cal M}_q(\mu)$ provides a rational parametrization of the polynomial Prony variety. It is a rational mapping of degree $2d -q -1.$

For the moment Prony curves $S^{\cal M}=S_{2d-2}^{\cal M}$, which are the straight lines in ${\cal M}_d$ parallel to the coordinate axis $Om_{2d-2}$, this restriction is linear in the last moment $m_{2d-1}$, and it is provided by the expression

\be\label{eq:param.of.Prony.curve}
c_i=C_1^i(\tilde \mu)m_{2d-1}+C_2^i(\tilde \mu),
\ee
where $\tilde \mu=(m_0,\ldots,m_{2d-2})$, and $C_1^i(\tilde \mu)$ and $C_2^i(\tilde \mu)$ are constant along the moment Prony curves $S^{\cal M}$.
\ep
An important fact is that the moment Hankel matrix $M_d(\mu)=M_d(\tilde \mu)$ {\it is constant along the moment Prony curves $S^{\cal M}(\mu)$}.

\section{Solvability over the reals}\label{Sec:Solv.Real.Hyperb}
\setcounter{equation}{0}



The requirement for all the amplitudes $A$ and the nodes $X$ of the reconstructed signal $F$ to be real is equivalent to requiring that all the moments
$\mu=\big(m_0(F),\allowbreak \ldots,m_{2d-1}(F) \big)\in {\cal M}_d$ are real, and that all the roots of the reconstructed polynomial $Q$ are real, i.e $Q$ is hyperbolic. As above,
we denote by $H_d\subset {\cal V}_d$ the set of hyperbolic polynomials.

We define the ``moment hyperbolicity set'' $\tilde H_d \subset {\cal M}_d$ as the set of all $\mu \in {\cal M}_d$ for which the Hankel image $HM(\mu)$ belongs to the hyperbolicity set $H_d\subset {\cal V}_d$. Equivalently,
$$
\tilde H_d = HM^{-1}(H_d).
$$
The following result is a partial case of the conditions of Prony solvability over the reals, obtained in \cite{goldman2018prony}:

\bt\label{th:solvsb.cond.real} (\cite{goldman2018prony}).
For a real moments vector $\mu \in {\cal M}_d$, with $\det M_d(\mu)$ nonzero, Prony system (\ref{eq:Moments.Prony11}) is solvable over the reals if and only if $\mu$ belongs to the moment hyperbolicity set $\tilde H_d \subset {\cal M}_d$.


\et

\subsection{Some statistical estimations for Prony solutions}\label{Prob.setting}

For a real signal $F$, if its moments vector $\mu$ was corrupted by the noise to $\mu'$, some roots of the reconstructed polynomial $Q=HM(\mu')$ could become complex. This makes the corresponding solution $F'$ unfeasible.

This situation is common in practice, and usually the complex roots of $Q$ are just projected to the real line. (In fact, in most of publications instead of real roots, the roots on the unit circle in the complex plane ${\mathbb C}$ are considered).

The same problem arises with the additional a priori known constraints on the feasible solutions $F$. (In particular, in most of applications we have a priori upper bounds on the nodes and amplitudes). We will denote by $Z\subset {\cal M}_d$  the set consisting of the moments of all the feasible signals $F$.

\smallskip

One of the most common {\it statistical estimations methods} is the maximum likelihood one (see e.g. \cite{plonka2018Optimal} and references therein).
Consider, for example, a Gaussian noise model $\mu' \sim {\cal N}(\mu,\Sigma)$ where
$\mu$ is unknown. Then the maximum likelihood estimator $\hat{\mu}(\mu')$ of $\mu$ is any point $z\in Z \subset {\cal M}_d$ that is nearest to the measurement $\mu'$.

\medskip

In Bayesian estimation, besides the assumed probability distribution for the noise (e.g. Gaussian), we also assume a prior probability distribution
of the moments vectors  (or of the feasible signals) with support on $Z$, and a fixed loss function $L(\hat\mu(\mu'),\mu)$.
Here the optimal Bayes estimator $\hat{\mu}(\mu')$ is given by the minimizer of the posterior risk
$$\hat{\mu}(\mu') = \inf_{\hat{\mu}\in Z} E[L(\hat{\mu},\mu)|\mu'] = \inf_{\hat{\mu} \in Z}\int_{Z} L(\hat\mu,\mu) f_{\mu|\mu'}(\mu)d\mu,$$
where $f_{\mu|\mu'}(\mu)$ is the conditional density of $\mu$ given the measurement $\mu'$.

\smallskip

Notice that minimisation is performed on an a priori known (and usually semi-algebraic) set $Z$. In our initial example $Z$ is the hyperbolicity domain $\tilde H_d\subset {\cal M}_d$. The study of such minimization problems is in the mainstream of Singularity Theory. Specifically, a rich geometric information on the hyperbolicity domain, available today, may be useful (see \cite{arnol1986hyperbolic,kostov1989geometric,kostov2011topics} and references therein). Another highly relevant topic in Singularity Theory is the study of singularities of maximal functions, cut loci, and related objects. Some ``old'' results are in \cite{Bryzgalova1977maximum,matov1982topological,yomdin1981local,yomdin1986functions,yomdin1983functions,zakalyukin1977singularities}\footnote{Let us notice that the proof of one of the main results in \cite{yomdin1981local} was incorrect, so the question remained open. Very recently a partial confirmation of the claim of \cite{yomdin1981local} was obtained in \cite{birbrair2017medial}.},
and in references therein. Some recent results are in \cite{reeve2012singularities,Che2017Sing}.

\section{Error amplification and Prony curves}\label{Sec:Error.Amplif.Prony}
\setcounter{equation}{0}

In this section we give a survey of recent results of \cite{akinshin2017error},
describing the geometry of error amplification in the case where the nodes of a signal $F$ form a cluster of size $h\ll 1$.
The central notion here is that of the {\it $\e$-error set $E_\e(F)$}.
\bd\label{def:error.set}
The error set $E_\e(F)\subset {\cal P}_d$ is the set consisting of all the signals $F' \in {\cal P}_d$ with
$$
|m_k(F')-m_k(F)|\le \e, \ k=0,\ldots,2d-1.
$$
\ed

In other words, $E_\e(F)$ comprises all the signals $F' \in {\cal P}_d$ which can appear in reconstruction of $F$ from its moments $\mu=(m_0,\ldots,m_{2d-1}),$ each moment $m_k$ corrupted by noise bounded by $\e$.\footnote{ In contrast with Section \ref{Prob.setting}, in \cite{akinshin2017error} and here we make no probabilistic assumptions on the noise.}

The goal here is a detailed understanding of the geometry of the error set $E_\e(F)$,
in the various cases where the nodes of $F$ near-collide.

\subsection{The model space}\label{Sec:model.signal}

For $F\in {\cal P}_d$, we denote by $I_F=[x_1,x_d]$, the minimal interval in
$\mathbb R$ containing all the nodes $x_1,\ldots,x_d$. We put $h(F)=\frac{1}{2}(x_d-x_1)$ to be the half of the length of $I_F$,
and put $\kappa(F)=\frac{1}{2}(x_1+x_d)$ to be the central point of $I_F$.

\smallskip

In case that $h(F)\ll 1$, we say that the nodes of $F$ form a cluster of size $h$ or simply that $F$ forms an
$h$-cluster.

For such signals $F$, consider the following ``normalization'': shifting the
interval $I_F$ to have its center at the origin, and then rescaling $I_F$ to the interval $[-1,1]$.
For this purpose we consider, for each $\kappa \in {\mathbb R}$ and $h>0$ the transformation
\be\label{eq:Ptransform}
\Psi_{\kappa,h}:{\cal P}_d\to {\cal P}_d,
\ee
defined by $(A,X)\to (A,\bar X),$ with
$$
\bar X=(\bar{x}_1,\ldots,\bar{x}_d), \ \ \bar{x}_j=\frac{1}{h}\left (x_j-\kappa\right), \ j=1,\ldots,d.
$$
For a given signal $F$ we put $h=h(F), \ \kappa=\kappa(F)$ and call the signal $G=\Psi_{\kappa,h}(F)$ the model signal for $F$.
Clearly, $h(G)=1$ and $\kappa(G)=0$. Explicitly $G$ is written as
$$
G(x)=\sum_{j=1}^{d}a_{j}\delta\left(x-\bar{x}_{j}\right).
$$
With a certain misuse of notations, we will denote the space ${\cal P}_d$ containing the model signals $G$ by
$\bar {\cal P}_d,$ and call it ``the model space''.
For $F\in {\cal P}_d$ and $G=\Psi_{\kappa,h}(F)$, the moments of $G$

\be\label{model.moments}
\bar m_k(F)= m_k(G)=\sum_{j=1}^d a_j\bar{x}_j^k, \ k=0,1,\ldots
\ee
are called the model moments of $F$.

\smallskip

For a given $F\in {\cal P}_d$ with the model signal $G=\Psi_{\kappa,h}(F),$
we denote by $\bar E_\e(F)$ the ``normalized" error set:
$$\bar E_\e(F) = \Psi_{\kappa,h}(E_\e(F)).$$
The set $\bar E_\e(F)$ represents the error set $E_\e(F)$ of $F$ in the model space $\bar {\cal P}_d$.
Note that $\bar E_\e(F)$ is simply a translated and rescaled version of $E_\e(F)$.

\smallskip

The reason for mapping a general signal $F$ into the model space is that in the case of the nodes $X$ forming
a cluster of size $h\ll  1$, the moment coordinates centered
at $F$,
$$\big(m_0(F')-m_0(F),\ldots,\allowbreak  m_{2d-1}(F')-m_{2d-1}(F)\big),$$
turn out to be ``stretched'' in some directions, up to the order $(\frac{1}{h})^{2d-1}$.
{\it In contrast, in the model space $\bar {\cal P}_d$,
the coordinates system $$\big(m_0(G')-m_0(G),\ldots, m_{2d-1}(G')-m_{2d-1}(G)\big)$$
is bi-Lipschitz equivalent to the standard coordinates $(A,\bar X)$ of $\bar {\cal P}_d$,
for all signals $G$ with ``well-separated nodes'' (see Theorem \ref{thm:coord.moments} below).}

\smallskip

Throughout this section we will always use the maximum norm $||\cdot||$ on ${\cal M}_d$ and on ${\cal P}_d$ and
on the nodes and amplitudes subspaces, ${\cal P}^X_d$ and ${\cal P}^A_d$ respectively.
Explicitly:\\
For $\mu = (\mu_0,\ldots,\mu_{2d-1}),\;\mu'=(\mu'_0,\ldots,\mu'_{2d-1}) \in {\cal M}_d$
$$
	\;\;\;\;||\mu' - \mu||=\max_{k=0,1,\ldots,2d-1} |\mu'_k-\mu_k|.
$$
For $F=(A,X),F'=(A',X')\in {\cal P}_d$,
$$
	||F-F'||=\max (||A-A'||,||X-X'||).
$$

%
%

\subsection{Sketch of the results}

We show that if the nodes of $F$ form a cluster of size $h\ll 1$ and $\e$ is of order $h^{2d-1}$ or less then:

\noindent{\it The $\e$-error set $\bar{E}_\e(F)$ is a ``curvilinear parallelepiped'' $\Pi$, which closely
follows the shape of the appropriate Prony varieties passing through $G$. The width of $\Pi$ in
the direction of the model moment coordinate $m_k(G')-m_k(G)$ is of order $\e h^{-k}.$}

Define the worst case reconstruction error of $F$ as
	$$
		\rho(F,\e)= \max_{F'\in E_\e(F)}||F'-F||.
	$$
	In a similar way we define $\rho_A(F,\e)$ and $\rho_X(F,\e)$ as the worst case errors in reconstruction of the amplitudes and the nodes of $F$, respectively:
	\begin{align*}
		\rho_A(F,\e)&= \max_{F'=(A',X')\in E_\e(F)}||A'-A||,\\
		\rho_X(F,\e)&= \max_{F'=(A',X')\in E_\e(F)}||X'-X||.
	\end{align*}
{\it We show that the worst case reconstruction error of the amplitudes $A$ and the signal $F$, $\rho_A(F,\e)$ and $\rho(F,\e)$,
are of order $\e h^{-2d+1}$, and, the worst case reconstruction error of the nodes $X$ is of order $\e h^{-2d+2}$.}

\smallskip

\noindent The above is shown in the following three steps:
\begin{enumerate}
	\item
		First we normalize the signal $F$ into its model signal $G=\Psi_{\kappa,h}(F)$,
		and describe in Theorem \ref{thm.error.set.geometry} the effect of this normalization on the
		image of the error set in ${\cal M}_d$.
		This theorem provides a description of the error set 	in the ``moment coordinates'', which
		are not, in general, equivalent to the coordinates of the signal space,
		because of the discussed singularities of the Prony mapping.
	\item
		The second step is to use a ``Quantitative Inverse Function Theorem''
		in order to show that the moment coordinates are bi-Lipschitz equivalent to the standard coordinates in signal space,
		in a sufficiently large domain around $G$. To get accurate constants,
		we improve in \cite{akinshin2017error} some estimates of the norm of the inverse Jacobian $JPM$ of the Prony mapping,
		obtained in \cite{batenkov2013accuracy}.
	\item
		Finally, in order to get accurate bounds for the worst case error
		{\it separately in the amplitudes $A$, and in the nodes $X$} of the reconstructed signal $F$,
		we provide in \cite{akinshin2017error} accurate estimates of the norm of the inverse Jacobian $JPM$ composed with the
		projections into the amplitudes and the nodes subspaces, ${\cal P}_d^A$ and ${\cal P}_d^X$, of ${\cal P}_d$.
\end{enumerate}
%
%

\subsection{The error set in the model signal space}\label{Sec:error.set.mod}

For any $G \in \bar{{\cal P}}_d$ and $\e,\alpha>0$ we denote by $\Pi_{\e,\alpha}(G)$ the ``curvilinear parallelepiped'' consisting of all $G' \in \bar{\cal P}_d$ satisfying
$$ |m_k(G')-m_k(G)|\leq \e \alpha^{k}, \ k=0,\ldots,2d-1.$$	

\smallskip

Notice that the Prony variety $S^{\cal P}_q(G)$ passing through $G$ is defined by the equations $m_k(G')=m_k(G), \ k=0,\ldots,q$, and therefore, in the moments coordinates $m_k(G')$ the parallelepiped $\Pi_{\e,h}(G)$ is $\e h^{-q}$ close to the Prony variety $S^{\cal P}_q(G)$.

\bt\label{thm.error.set.geometry}
Let $F\in {\cal P}_d$ form a cluster of size $h=h(F)$ and let $\kappa=\kappa(F)$ be the center of the cluster. Let $G=\Psi_{\kappa,h}(F)$ be the model signal for
$F$. Set $\e'=(1+|\kappa|)^{-2d+1} \e$ and $h'=\frac{h}{1+|\kappa|}$. Then for any $\e>0$, the error set $\bar{E}_\e(F)$ is
bounded between the following two parallelepipeds:
  		$$
			 \Pi_{ \e',\frac{1}{h}}(G) \subset \bar E_\e(F) \subset \Pi_{\e,\frac{1}{h'}}(G).
		$$
Specifically, for $\kappa=\kappa(F)=0$,
		$$
			\bar E_\e(F)=\Pi_{\e,h}(G).
		$$
\et
Theorem \ref{thm.error.set.geometry} holds without any assumptions on the mutual relation of $\e$ and $h$, or on the distances between the nodes of $F$. It implies the following fact:
{\it the Prony varieties $S^{\cal P}_q(G)$ form a ``skeleton'' of the error set $\bar E_\e(F)$, and, in case when $\e$ and $h$ tend to zero at a certain rate, $S^{\cal P}_q(G)$ are
the limits of $\bar E_\e(F)$}.

\smallskip

\begin{figure}
		\centering
		\includegraphics[scale=0.70]{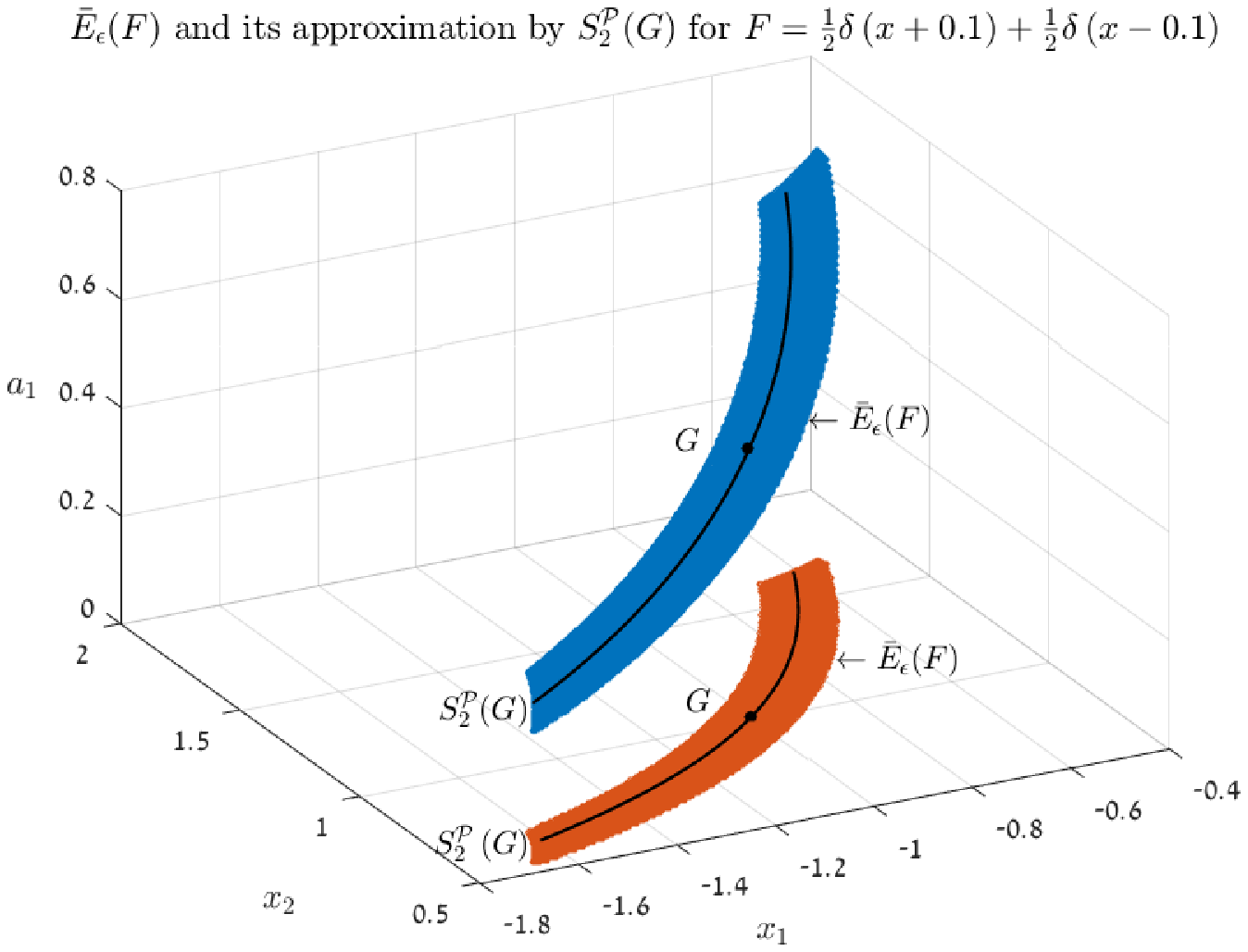}
		\caption{The projections of the error set $\bar{E}_{\e}(F)$ and a section of the Prony
		curve $S^{\cal P}_2(G)$, for $h=0.1$ and $\epsilon = h^3$.}
		\label{fig.h01}
		\includegraphics[scale=0.70]{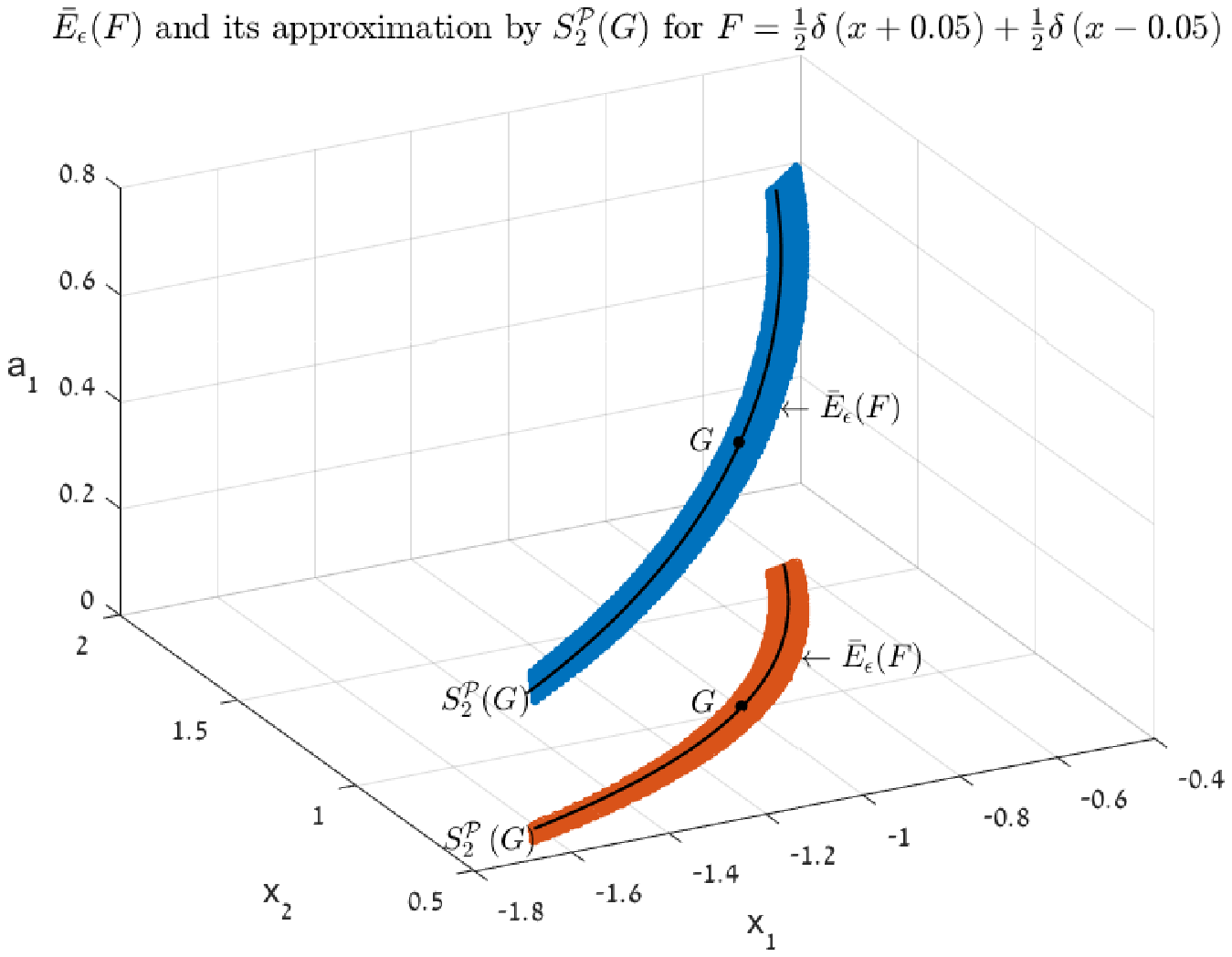}
		\caption{The error set $\bar{E}_{\e}(F)$ and a section of $S^{\cal P}_2(G)$ for
		$h=0.05$ and $\epsilon = h^3$. Note the convergence of
		$\bar{E}_{\e}(F)$ into $S^{\cal P}_2(G)$.}
		\label{fig.h05}
\end{figure}

Figures \ref{fig.h01} and \ref{fig.h05} illustrate the case $d=2, q=2d-2=2$ of Theorem \ref{thm.error.set.geometry}.

\subsection{Applying quantitative Inverse Function Theorem}\label{Sec:appl.QIFT}

In order to apply this theorem, we have to make explicit assumptions on the separation of the nodes $X$ of the signal $G$, and on the size of its amplitudes $A$:

Assume that the nodes $x_1,\ldots,x_d$ of a signal $G \in \bar{{\cal P}}_d$ belong to the interval $I=[-1,1]$,
and for a certain $\eta$ with $0<\eta\leq \frac{2}{d-1}$, $d>1$,
the distance between the neighboring nodes $x_j,x_{j+1}, \ j=1,\ldots,d-1,$ is at least $\eta$.
We also assume that for certain $m,M$ with $0<m<M$, the amplitudes $a_1,\ldots,a_d$ satisfy $m\leq |a_j|\leq M, \ j=1,\ldots,d$. We call such signals $(\eta,m,M)$-regular.

\smallskip

We want to show that for an $(\eta,m,M)$-regular signal $G\in \bar {\cal P}_d$ the moment coordinates
$m_0(G')-m_0(G),\ldots,m_{2d-1}(G')-m_{2d-1}(G)$ indeed form a coordinate system
near $G$, which agrees with the standard coordinates $A,\bar X$ on $\bar {\cal P}_d$.

\bd\label{def:moment.coord.dist}
For $G$ a regular signal as above, and $G'$ denoting signals near $G$,
the moment coordinates are the functions $f_k(G')=m_k(G')-m_k(G),\; k=0,...,2d-1$.
The moment metric $d(G',G'')$ on $\bar {\cal P}_d$ is defined through the moment coordinates as
$$
d(G',G'')=\max_{k=0}^{2d-1}|m_k(G'')-m_k(G')|.
$$
\ed

For any $\nu \in {\cal M}_d$ and $R >0$ denote by $Q_{R}(\nu) \subset {\cal M}_d$ the cube of radius $R$
\be\label{eq.epsilon.cube}
	Q_{R}(\nu)=\{\nu'=(\nu'_0, \ldots, \nu'_{2d-1}) \in {\cal M}_d, \ |\nu'_k - \nu_k|\le R, \ k=0,1,\ldots,2d-1\}.
\ee

\bt\label{thm:coord.moments}
Let $G\in \bar {\cal P}_d$ be an $(\eta,m,M)$ regular signal and $\nu = PM(G)$. Then there are constants $R,C_1,C_2$, depending only on $d,\eta,m,M$,
given explicitly in \cite{akinshin2017error}, such that:

\begin{enumerate}
	\item The inverse mapping $PM^{-1}$ exists on $Q_R(\nu)$ and provides a diffeomorphism of
	$Q_R(\nu)$ to $\Omega_R(G)= PM^{-1}(Q_R(\nu))$.
	\item The moment metric $d(G',G'')$ is bi-Lipschitz equivalent on $\O_R(G)$ to the maximum metric $||G''-G'||$ in $\bar{\cal{P}}_d$:
$$
	C_1 \ d(G',G'')\le ||G''-G'||\le C_2 \ d(G',G'').
$$
\end{enumerate}

%

\et

Assume now that the measurement error $\e \le R h'^{2d-1}$, with $h'=\frac{h}{1+|\kappa|}$ as
in Theorem \ref{thm.error.set.geometry}. Then
	$$PM(\bar{E}_\e(F)) \subseteq PM(\Pi_{\e,\frac{1}{h'}}(G)) \subset Q_R(PM(G)).$$
Combing Theorems \ref{thm.error.set.geometry} and \ref{thm:coord.moments}
we obtain that the error set $\bar{E}_\e(F)$ is a ``deformed'' paralelipiped in $\bar{{\cal P}}_d$ as
illustrated in figures \ref{fig.h01} and \ref{fig.h05} above.

\smallskip

We use regular signals $G$ as above, to model signals with a
``regular cluster'': For $F\in {\cal P}_d$ with $h=h(F)$ and $\kappa=\kappa(F)$, we say that $F$ forms an $(h,\kappa,\eta,m,M)$-regular cluster if $G=\Psi_{\kappa,h}(F)$ is an $(\eta,m,M)$-regular signal.

The next theorem shows that the $\e$-error set is tightly concentrated around the Prony varieties.

\bd
	For each $0 \le q \le 2d-1$ denote by $S^{{\cal P}}_{q,\e,\alpha}(G)$ the part of the Prony variety $S^{{\cal P}}_q(G),$
	consisting of all signals $G'\in S^{{\cal P}}_q(G)$ with
	$$
		|m_k(G')-m_k(G)|\leq \e \alpha^{k}, \ k=q+1,\ldots,2d-1.
	$$
\ed
%
%

\bt\label{thm:distance.to.Sq}
Let $F\in {\cal P}_d$ form an $(h,\kappa,\eta,m,M)$-regular cluster and let $G=\Psi_{\kappa,h}(F)$ be the model signal for
$F$. Set $h'=\frac{h}{1+|\kappa|}$. Then
for any $\e \le R h'^{2d-1}$, the error set $\bar{E}_{\e}(F)$
is contained within the $\Delta_q$-neighborhood of the part of the Prony variety
$S^{{\cal P}}_{q,\e,\frac{1}{h'}}(G)$, for
   		$$
 			\Delta_q=C_2 \left(\frac{1}{h'}\right)^{q}\e.
 		$$
 		The constants $R,C_2$ are defined in Theorem \ref{thm:coord.moments} above.
\et


\subsection{Worst case reconstruction error}

We now present lower and upper bounds, of the same order, for the worst case reconstruction error
$\rho(F,\e),$ defined, as above, by:
$$
\rho(F,\e)= \max_{F'\in E_\e(F)}||F'-F||.
$$
We state separate bounds for $\rho_A(F,\e)$ and $\rho_X(F,\e)$ - the worst case errors in reconstruction of the amplitudes $A=(a_1,\ldots,a_d)$ and of the nodes $X=(x_1,\dots,x_d)$ of $F$:
$$
\rho_A(F,\e)= \max_{F'\in E_\e(F)}||A'-A||, \ \rho_X(F,\e)= \max_{F'\in E_\e(F)}||X'-X||.
$$
\bt\label{thm.upper.bound}[Upper bound]
Let $F\in {\cal P}_d$ form an $(h,\kappa,\eta,m,M)$-regular cluster.
Then for each positive $\e\le \left(\frac{h}{1+|\kappa|}\right)^{2d-1}R$ the following bounds for the worst case reconstruction errors are
valid:
$$
	 \rho(F,\e),\; \rho_A(F,\e)  \le C_2 \left(\frac{1+|\kappa|}{h}\right)^{2d-1} \e,\tab \rho_X(F,\e)  \le C_2 h \left(\frac{1+|\kappa|}{h}\right)^{2d-1} \e,
$$
where $C_2, R$ are the constants defined in Theorem \ref{thm:coord.moments}.
\et

\bt\label{thm.lower.bound}[Lower bound]
Let $F\in {\cal P}$ form an $(h,\kappa,\eta,m,M)$-regular cluster then:
\begin{enumerate}
  \item For each positive $\e\le C_3 h^{2d-1}$ we have the following lower bound on the worst case reconstruction error of the nodes of $F$
  		$$K_1 \e \left(\frac {1}{h}\right)^{2d-2}\le \rho_X(F,\e).$$
  \item For each positive $\e\le C_4 h^{2d-1}$ we have the following lower bound on the worst case reconstruction error of $F$ and of the amplitudes of $F$
     $$K_2 \e \left(\frac{1}{h}\right)^{2d-1}\le \rho(F,\e),\;\rho_A(F,\e).$$
\end{enumerate}
	Above, $K_1, K_2, C_3, C_4$ are constants not depending on $h$
	given explicitly in \cite{akinshin2017error}.
\et

The lower and upper bounds given above are a special case of a more general result.
In \cite{akinshin2017error} (Theorem 4.4) it is shown that the Prony variety $S^{\cal P}_{q}(G)$ can be reconstructed
from the moment measurements $\mu' \in{\cal M}_d$ with improved accuracy of order $\e h^{-q}$.

\section{Prony Scenarios}\label{Sec:Prony.Scenarios}
\setcounter{equation}{0}


We keep the assumption that the nodes of our signal $F$ form a regular cluster of a size $h\ll 1$. By Theorem \ref{thm:distance.to.Sq}, the signal Prony curve $S^{\cal P}(\mu)$ approximates the error set $E_\e(F)$ with the accuracy of order $\e h^{-2d+2}$. Note that the accuracy of point solution is of order $\e h^{-2d+1}$.
Thus, the Prony curve $S^{\cal P}(\mu)$  provides a rather accurate prediction of the possible behavior of all the noisy reconstructions of $F$.

In an actual solution procedure, the ``true'' Prony curve $S^{\cal P}(\mu)$ is not known. But from the noisy measurements $\mu'=(m'_0,\ldots,m'_{2d-1})$ we can reconstruct the Prony curve $S^{\cal P}(\mu')$.
This curve, by Theorem 4.4 in \cite{akinshin2017error}, approximates the ``true curve'' $S^{\cal P}(\mu)$ with an accuracy of the same (improved) order of $\e h^{-2d+2}$, with which $S^{\cal P}(\mu)$ approximates $E_\e(F)$. Therefore, we can consider this known curve $S^{\cal P}(\mu')$ as a prediction (or a ``scenario'') for all the noisy reconstructions of $F$.

\smallskip

{\it Moreover, if we neglect possible errors of order $\e h^{-2d+2}$, we can restrict the search of the optimal Prony solution (by any method, in particular, via statistical estimations) to the curve $S^{\cal P}(\mu')$.}


\smallskip



We do not try to give here a rigorous definition of the ``Prony scenario''. Informally, this is a collection of data on the Prony curve $S^{\cal P}(\mu')$, which is necessary in order to find the optimal Prony solution on this curve, taking into account the available a priori constraints. Certainly we need an accurate description of the behavior of the nodes $x_j$ and the amplitudes $a_j$ along $S^{\cal P}(\mu')$ (or, better, along the polynomial Prony curve $S^{\cal V}(\mu')\subset {\cal V}_d$), including description of the intersection of $S^{\cal V}(\mu')$ with the hyperbolicity set $H_d$.


Some general results in this direction were obtained in \cite{goldman2018prony}:

\bt\label{Thm:amplitudes.yoni} (\cite{goldman2018prony})
Assume that the matrix $M_d(\mu')$ is non-degenerate. Then in each case where the nodes $x_i,x_j$ collide on $S^{\cal P}(\mu')$, the amplitudes $a_i$ and $a_j$ tend to infinity.
\et
\bt\label{Thm:nodes.yoni} (\cite{goldman2018prony})
Assume that the matrix $M_d(\mu')$ is non-degenerate, as well as its upper-left $(d-1)\times(d-1)$ minor. Then on each unbounded component of $S^{{\cal P},X}(\mu')$, for the coordinate $m_{2d-1}$ on $S^{\cal P}(\mu')$ tending to infinity, exactly one node ($x_1$ or $x_d$) tends to infinity, while the rest of the nodes remain bounded.
\et
The polynomial Prony curves $S^{\cal V}$ can be considered as polynomials pencils. Some important results on the behavior of the real roots in polynomial pencils are provided in \cite{borcea2004classifying,kurdyka201nuij}.
The result of \cite{parusinski2013regularity} describing the behavior of roots in smooth 1-parametric families of polynomials may also be relevant. 
These results naturally enter the framework of the Prony scenarios, and in \cite{goldman2018prony} we provide their more detailed treatment.

\smallskip





\section{Some open questions}\label{Sec:Open.Quest}
\setcounter{equation}{0}

We would like to specify some open problems in the line of this paper. Mostly they concern the structure of the Prony varieties in the areas not covered by the inverse function theorem (Theorem \ref{thm:coord.moments} above).


\smallskip

\noindent {\it 1. Description of the global topology and geometry of the Prony varieties}. In the topological study of the Vandermonde varieties in \cite{arnol1986hyperbolic,kostov1989geometric} certain natural
Morse functions were used. Can this method be extended to the Prony varieties?

\smallskip

On the other hand, an explicit parametrization of the Prony varieties, described in Section \ref{Sec:Prony.Var.in.Three.Spaces} above, reduces the problem to the study of the intersections of the affine subspaces in the polynomial space with the hyperbolic set $H_d$ (which is motivated also by the considerations in Section \ref{Sec:Prony.Scenarios} above). This study looks natural also from the point of view of Singularity Theory.

\smallskip

\noindent {\it 2. Understanding connections between the Prony and the Vandermonde varieties}. The last are the fibers of a natural projection of the corresponding Prony varieties
to the amplitudes. Is this projection regular? What topological information on the Prony varieties can be obtained from the known properties of the Vandermonde ones? Can
information available on the Prony varieties (in particular, their explicit parametrization, see Section \ref{Sec:Prony.Var.in.Three.Spaces} above) be useful in study of the
Vandermonde ones?
 
\smallskip

\noindent {\it 3. Behavior of the nodes $x_1,\ldots,x_d$ on the Prony varieties $S^{\cal P}_q(\mu) \subset {\cal P}_d$  near the collision singularities}. It would be important to describe an accurate asymptotic behavior of the distances between the colliding nodes as we approach the collision point. This question can be split into two: investigation of the intersection of the affine varieties $S^{\cal V}_q(\mu)\subset {\cal V}_d$ with the boundary of the hyperbolic set $H_d$, and investigation of the behavior of the root mapping $RM$ near the boundary of $H_d$.

\smallskip

Already the case of the Prony curves is important and non-trivial.



\smallskip

\noindent {\it 4. Behavior of the amplitudes $a_1,\ldots,a_d$ on the Prony varieties $S^{\cal P}_q(\mu) \subset {\cal P}_d$ near the collision singularities}. 
In the case of the Prony curve, i.e. $q=2d-2$, Theorem \ref{Thm:amplitudes.yoni} above gives conditions under which these amplitudes necessarily tend to infinity.
It would be important 
to describe the accurate asymptotic behavior of the amplitudes as we approach the collision point. We expect that this question can be treated via methods from the classical Moment theory, combined with the techniques of ``bases of finite differences'' developed in \cite{batenkov2013geometry,yomdin2010singularities}. Also here the case of the Prony curves is important.

\smallskip



\noindent {\it 5. Extending the description of the Prony varieties, and of the error amplification patterns, to multi-cluster nodes configurations}.
This is a natural setting in robust inversion of the Prony system. 
Generalized Prony methods as well as other reconstruction methods typically reduce each cluster to a single node, 
thus forming a ``reduced Prony system''. It is important to estimate the accuracy of such an approximation (see \cite{goldman2018algebraic} for some steps in this direction).

\smallskip

Because of the role of the Prony varieties in the analysis of the error amplification patterns, a natural question is: {\it To what extent do the Prony
varieties of the reduced Prony system approximate the varieties of the ``true'' multi-cluster system?}

\bibliographystyle{myplain}
\bibliography{bib}{}

\end{document}